%% file: main.tex
\documentclass[reqno, 12pt]{amsart}
\usepackage[margin=1in]{geometry}
\pdfoutput=1
\usepackage[T2A]{fontenc}
\usepackage[english]{babel}

\usepackage{amssymb}
\usepackage{amsthm}
\usepackage{amsfonts}
\usepackage{amsmath}
\usepackage{amscd}
\usepackage{bbm}
\usepackage{eucal}
\usepackage{faktor}
\usepackage{xfrac}
\usepackage{braket}
\usepackage{relsize}
\usepackage{nccmath}
\usepackage{mathtools}
\usepackage{comment}
\usepackage{verbatim}
\usepackage{wasysym}
\usepackage{graphicx}
\usepackage{import}
\usepackage{float}
\usepackage{mathrsfs}
\usepackage{color}
\usepackage{xparse}
\usepackage{chngcntr}
\usepackage{mwe}
\usepackage{enumitem}
\usepackage{tikz}
\usepackage{graphicx}
\usepackage{wrapfig}
\usepackage{lipsum}
\usepackage{dsfont}
\usepackage{etoolbox}
\usepackage{stackrel}
\usepackage{bm}
\usepackage{upgreek}
\usepackage{url}
\usepackage{afterpage}
\usepackage{scalerel}[2016/12/29]
\usepackage{genyoungtabtikz}

\usepackage[
backref=true,
backrefstyle=none,
backend=biber,
style=alphabetic,
citestyle=alphabetic
]{biblatex}

\DefineBibliographyStrings{english}{
  backrefpage = {cited on page},
  backrefpages = {cited on pages},
}

\usepackage[figurename=Figure.]{caption}

\usetikzlibrary{patterns}
\usetikzlibrary{shapes}

\usepackage[linktocpage=true]{hyperref}

\mathtoolsset{showonlyrefs,showmanualtags}

\allowdisplaybreaks

\newcommand{\quot}[2]{{\raisebox{.3em}{$#1$}\bigl/\raisebox{-.3em}{$\!#2$}}}

\newcommand\T{\mathcal{T}}
\newcommand\CC{\mathbb{C}}
\newcommand\Pas{\mathbb{P}}
\newcommand\NN{\mathbb{N}}
\newcommand\YY{\mathbb{Y}}

\newcommand\Ga{\Gamma}

\newcommand\ka{\varkappa}

\newcommand\la{\lambda}

\newcommand\wt{\widetilde}
\newcommand\wh{\widehat}

\newcommand\summ{\sum\limits}
\newcommand\limm{\lim\limits}

\makeatletter
\newcommand*{\isomorphism}{%
  \mathrel{%
    \mathpalette\@isomorphism{}%
  }%
}
\newcommand*{\@isomorphism}[2]{%
  \sbox0{$#1\simeq$}%
  \sbox2{$#1\sim$}%
  \dimen@=\ht0 %
  \advance\dimen@ by -\ht2 %
  \sbox0{%
    \lower2.5\dimen@\hbox{%
      $\m@th#1\relbar\isomorphism@joinrel\rightarrow$%
    }%
  }%
  \rlap{%
    \hbox to \wd0{%
      \hfill\raise\dimen@\hbox{$\m@th#1\sim$}\hfill
    }%
  }%
  \copy0 %
}
\newcommand*{\isomorphism@joinrel}{%
  \mathrel{%
    \mkern-3.4mu %
    \mkern-1mu %
    \nonscript\mkern1mu %
  }%
}
\makeatother

\newcommand\const{\operatorname{const}}
\newcommand\supp{\operatorname{supp}}
\newcommand{\spn}[2]{\operatorname{span}_{#1}\left(#2\right)}

\renewcommand\H{\mathcal{H}}
\DeclareMathOperator{\FH}{\mathcal{FH}}
\DeclareMathOperator{\exH}{\H_{ex}}
\DeclareMathOperator{\exFH}{\FH_{ex}}
\DeclareMathOperator{\exHo}{\H_{ex}^{\circ}}
\DeclareMathOperator{\exFHo}{\FH_{ex}^{\circ}}
\newcommand\sym{\mathit{Sym}}

\DeclareMathOperator{\K}{K_0}
\DeclareMathOperator{\KK}{K_0^+}

\renewcommand{\binom}[2]{\begin{pmatrix}
#1\\ 
#2
\end{pmatrix}
}
\newtheorem{theorem}{Theorem}
\newtheorem{proposition}[theorem]{Proposition}
\newtheorem{corollary}[theorem]{Corollary}
\newtheorem{lemma}[theorem]{Lemma}

\newtheorem*{notation}{Notation}

\theoremstyle{definition}
\newtheorem{definition}[theorem]{Definition}
\newtheorem{definition-proposition}{Definition-Proposition}
\newtheorem{remark}[theorem]{Remark}

\newtheorem{notation-definition}[theorem]{Notation-Definition}

\numberwithin{definition-proposition}{section}
\numberwithin{theorem}{section}

\title{Semifinite harmonic functions on the direct product of graded graphs}

\author{Pavel Nikitin}

\address{
Laboratory of Representation Theory and Dynamical Systems,
St. Petersburg Department of Steklov Mathematical Institute
of Russian Academy of Sciences
191023, Fontanka 27, St. Petersburg, Russia.
}

\email{pnikitin0103@yahoo.co.uk}

\author{Nikita Safonkin}

\address{Skolkovo Institute of Science and Technology, Moscow, Russia \& National Research University Higher School of Economics, Moscow, Russia.}

\email{safonkin.nik@gmail.com}

\begin{document}

\begin{abstract}
    Indecomposible semifinite harmonic functions on the direct product of graded graphs are classified. As a particular case, the full list of indecomposible traces for the infinite inverse symmetric semigroup is obtained.  
\end{abstract}

    \maketitle
    \setcounter{tocdepth}{2}
    \tableofcontents

\import{}{Notes}

    \printbibliography
\end{document}

%% file: Notes.tex
\section{Introduction}

In the present paper we describe real-valued nonnegative indecomposable semifinite harmonic functions on the direct product of graded graphs. Infinite values are allowed for semifinite harmonic functions, but should be approximated by finite ones, see Definition \ref{def5}. Finite harmonic functions on a graded graph represent $II_1$--factor representations of the corresponding approximately finite-dimensional algebra, and are known in many interesting cases. 
And \emph{semifinite} harmonic functions represent a subclass of $II_\infty$--factor representations where a classification should be still possible for most cases. 

Currently the main tool to study semifinite traces is Wassermann's method~\cite[III \S6]{wassermann1981}, \cite{boyer83}; see also \cites{safonkin21,safon22}. But in the initial paper by Kerov and Vershik~\cite{vershik_kerov1980} the use of the \emph{ergodic method} was proposed. For the direct product of graded graphs we have the product structure on the space of paths, and it is reasonable to expect that the ergodic method works nicely in this case. Nevertheless, it turns out that it is more convenient to use Wasserman's method for such graphs as well. Namely, we restrict out attention to the finiteness ideal, where the classification is known, and then use a bijection between the indecomposable harmonic functions and their restrictions to an ideal to obtain the main result.

We denote the set of all indecomposable finite and semifinite harmonic functions on a graded graph $\Ga$ by $\exH(\Ga)$. 

\begin{theorem}[Main theorem]\label{th1}
    Let $\Ga_1$ and $\Ga_2$ be graded graphs and $\varphi\in\exH(\Ga_1\times\Ga_2)$, then only one of the following situations can occur: 
    \begin{enumerate}[label=\theenumi)]
    \item\label{th1case1} There exist $\varphi_1\in\exH(\Ga_1)$, $\varphi_2\in\exH(\Ga_2)$ and real positive numbers $w_1,w_2$ with $w_1+w_2=1$ such that  
\begin{equation}\label{f09}
\varphi(\la,\mu)=w_1^{|\la|}w_2^{|\mu|}\varphi_1(\la)\varphi_2(\mu).   
\end{equation} 
Moreover, these $\varphi_1$ and $\varphi_2$ are defined uniquely up to multiplicative constants.\newline
We agree that $0\cdot (+\infty)=0$.

    \item\label{th1case2} There exist $\varphi_1\in\exH(\Ga_1)$ and $\nu_2\in\Ga_2$ such that
    \begin{equation}\label{f10}
        \varphi(\la,\mu)=\begin{cases}
            0,&\text{if}\ \mu\not\leq\nu_2,\\
            +\infty,&\text{if}\ \mu<\nu_2,\\
            \varphi_1(\la),&\text{if}\ \mu=\nu_2.
        \end{cases}    \end{equation}

    \item\label{th1case3} There exist $\nu_1\in\Ga_1$ and $\varphi_2\in\exH(\Ga_2)$ such that 
    \begin{equation}\label{f11}
        \varphi(\la,\mu)=\begin{cases}
            0,&\text{if}\ \la\not\leq\nu_1,\\
            +\infty,&\text{if}\ \la<\nu_1,\\
            \varphi_2(\mu),&\text{if}\ \la=\nu_1.
        \end{cases}    \end{equation}
    \end{enumerate} 
    Furthermore, every harmonic function on $\Ga_1\times\Ga_2$ of the form \ref{th1case1}, \ref{th1case2}, or \ref{th1case3} is finite or semifinite, and indecomposable.
\end{theorem}

Direct products of graded graphs arise naturally in the study of semifinite traces, see~\cite[Lemma 4.15]{safon22}, as well as in the representation theory of inverse symmetric semigroups~\cite{Vershik-Nikitin, Halverson}. For direct products finite traces have been classified \cite{Vershik-Nikitin,safonkin21}, see also Theorem~\ref{th7}, and one of the main motivations for the present paper was to obtain a similar picture in the semifinite case. For the inverse symmetric semigroups the graph in question is the product $\NN \times \YY$, where $\YY$ is the Young graph and $\NN$ is the representation of the natural numbers as a half-line Bratteli diagram. The main theorem combined with the known results for the Young graph, gives the full characterization of the indecomposable semifinite traces for the infinite symmetric inverse semigroup, see  Proposition~\ref{prop:inv-semigroup_semifin-traces}.

The paper is organized as follows. In Section \ref{sec2} we recall the main notions related to branching graphs and semifinite harmonic functions on them. In Section \ref{sec3} we define the direct product of branching graphs and prove our main result, Theorem \ref{th1}. In Section \ref{sec4} we discuss the case of the infinite inverse symmetric semigroup. In Appendix \ref{appendix} we discuss the proof of an analog of the main result for the case of finite harmonic functions and some related facts.

\textbf{Ackowledgements}

This work is supported by the Russian Science Foundation under grant No.~21-11-00152. The second author was partially supported by the Basic Research Program at the HSE University. We thank Anatoly Vershik and Grigori Olshanski for useful comments and important discussions. We are greatful to the Leonhard Euler International Mathematical Institute at Saint Petersburg for the hospitality in December 2021 and a great opportunity for collaboration.

\section{Preliminaries on graded graphs}\label{sec2}
In this section we recall main notions related to  graded graphs and semifinite harmonic functions on them.  A detailed combinatorial exposition can be found in \cite{safonkin21}, see also original papers on traces on AF-algebras \cite{strat_voic1975,wassermann1981,Bratteli1972} and \cite{vershik_kerov83,versh_ker_85,kerov_vershik1990}.

\begin{definition}\label{def2}
    By a \textit{graded graph} we mean a pair $\left(\Ga,\ka\right)$, where $\Ga$ is a graded set $\Ga=\bigsqcup\limits_{n\geq 0}\Ga_n$, $\Ga_n$ are finite sets and $\ka$ is a function $\Ga\times \Ga\rightarrow \mathbb{R}_{\geq0}$, that satisfies the following constraints:
    \begin{enumerate}[label=\theenumi)]
    \item\label{def2cond1} if $\la\in\Ga_n$ and $\mu\in\Ga_m$, then $\ka(\la,\mu)=0$ for $m-n\neq 1$. 
    
    \item\label{def2cond2} for any vertex $\la\in\Ga_n$ there exists $\mu\in\Ga_{n+1}$ with $\ka(\la,\mu)\neq 0$.
\end{enumerate}
Edges of the graded graph $(\Ga,\ka)$ are, by definition, pairs of vertices $(\la,\mu)$ with $\ka(\la,\mu)>0$. Then we may treat $\ka(\la,\mu)$ as a formal multiplicity of the edge.
\end{definition}

If $\la\in\Ga_n$, then the number $n$ is uniquely defined. We denote it by $|\la|$. We write $\la\nearrow \mu,$ if $|\mu|-|\la|=1$ and $\ka(\la,\mu)\neq 0$. In this case we say that \textit{there is an edge from $\la$ to $\mu$ of multiplicity $\ka(\la,\mu)$}.

Condition \ref{def2cond1} from Definition \ref{def2} means that we allow edges only between adjacent levels and condition \ref{def2cond2} means that each vertex must be connected by an edge with some vertex from the higher level.

\textit{A path} in a graded graph $\Ga$ is a (finite or infinite) sequence of vertices  $\la_1,\la_2,\la_3,\ldots$ such that $\la_i\nearrow \la_{i+1}$ for every $i$. We will write $\nu>\mu$ if $|\nu|>|\mu|$ and there is a path that connects $\mu$ and $\nu$. We write $\nu\geq \mu$, if $\nu=\mu$ or $\nu>\mu$. Relation $\geq$ turns $\Ga$ into a poset.

Let $\mu,\nu\in\Ga$ and $|\nu|-|\mu|=n\geq 1$. Then the following expression
\begin{equation}
\dim(\mu,\nu)=\summ_{\substack{\la_0,\ldots,\la_n\in\Ga:\\ \mu=\la_0\nearrow\la_1\nearrow\ldots\nearrow\la_{n-1}\nearrow\la_n=\nu}}\ka(\la_0,\la_1)\ka(\la_1,\la_2)\ldots\ka(\la_{n-1},\la_n).
\end{equation}

is the "weighted" number of paths from $\mu$ to $\nu$. By definition we also set $\dim(\mu,\mu)=1$ and $\dim(\mu,\nu)=0$, if $\nu\not\geq\mu$. The function $\dim(\cdot,\cdot)\colon \Ga\times\Ga\rightarrow \mathbb{R}_{\geq 0}$ is called \textit{the shifted dimension}. 

\begin{definition}
    A \emph{branching graph} (or a \emph{Bratteli diagram}) is defined as a graded graph $(\Ga,\ka)$ that satisfies the following conditions
    \begin{itemize}
        \item $\Ga_0=\{\diameter\}$ is a singleton,
        
        \item for any $\la\in\Ga_n$ with $n\geq 1$ there exists $\mu\in\Ga_{n-1}$ such that $\mu\nearrow\la$.
    \end{itemize}
\end{definition}

For a branching graph $(\Ga,\ka)$ we denote the expression $\dim(\diameter,\la)$ by $\dim(\la)$ and call it \textit{the dimension of} $\la$.

\subsection{Semifinite harmonic functions}

Let $(\Ga,\ka)$ be a graded graph.

\begin{definition}
 A function $\varphi\colon \Ga\rightarrow \mathbb{R}_{\geq 0}\cup\{+\infty\}$ is called \emph{harmonic}, if it satisfies the following condition 
 $$\varphi(\la)=\summ_{\mu:\la\nearrow \mu}\ka(\la,\mu)\varphi(\mu).$$ 
 The set of all vertices $\la\in\Ga$ with $\varphi(\la)<+\infty$ is called the \textit{finiteness ideal} of $\varphi$. We denote the zero set $\left\{\la\in\Ga\; \middle|\; \varphi(\la)=0\right\}$ by $\ker{\varphi}$, and call it the \textit{kernel} or \textit{zero ideal}, and the \textit{support} $\{\la\in\Ga\mid \varphi(\la)>0\}$ by $\supp{\varphi}$.
\end{definition}

Note that there is an obvious bijection between harmonic functions on a graded graph $\Ga$ with the given support $J$ and strictly positive harmonic functions on $J$. 

The symbol $\K(\Ga)$ stands for the $\mathbb{R}$-vector space spanned by the vertices of $\Ga$ subject to the following relations $$\la=\summ_{\mu:\la\nearrow \mu}\ka(\la,\mu)\cdot \mu,\ \ \forall \la\in\Ga.$$ 
The symbol $\KK(\Ga)$ denotes the positive cone in $\K(\Ga)$, generated by the vertices of $\Ga$, i.e. $\KK(\Ga)=\spn{\mathbb{R}_{\geq0}}{\la\mid\la\in\Ga}$. The partial order, defined by the cone $\KK(\Ga)$, is denoted by $\geq_K$. That is $a\geq_K b\iff a-b\in \KK(\Ga)$. For instance, if $\la\geq \mu$, then $\mu\geq_K\dim(\mu,\la)\cdot\la$.

The $\mathbb{R}_{\geq0}$-linear map $\KK(\Ga)\rightarrow \mathbb{R}_{\geq0}\cup\{+\infty\}$, defined by a harmonic function $\varphi$, will be denoted by the same letter $\varphi$. 

\begin{definition}\label{def5}
A harmonic function $\varphi$ is called \emph{semifinite}, if it is not finite and for any $a\in\KK(\Ga)$ the map $\varphi\colon\KK(\Ga)\rightarrow \mathbb{R}_{\geq0}\cup\{+\infty\}$ enjoys the following property 
\begin{equation}\label{f28}
    \varphi(a)=\sup\limits_{\substack{b\in\KK(\Ga)\colon  b\leq_K a,\\ \varphi(b)<+\infty}}\varphi(b).
\end{equation}
\end{definition}

\begin{definition}
A semifinite harmonic function $\varphi$ is called \textit{indecomposable}, if for any finite or semifinite harmonic function $\varphi'$ which does not vanish identically on the finiteness ideal of $\varphi$ and satisfies the inequality $\varphi'\leq \varphi$ we have $\varphi'=\const\cdot\varphi$ on the finiteness ideal of $\varphi$.
\end{definition}

\begin{definition}
    A subset of vertices $I$ of a graded graph $\Ga$ is called an \textit{ideal}, if for any vertices $\la\in I$ and $\mu\in\Ga$ such that $\mu>\la$ we have $\mu\in I$. A subset $J\subset \Ga$ is called a \textit{coideal}, if for any vertices $\la\in J$ and $\mu\in\Ga$ such that $\mu<\la$ we have $\mu\in J$.
\end{definition}

Note that the kernel of a harmonic function is an ideal and the support is a coideal.

\begin{definition}\label{def4}
A graded graph $\Ga$ is called \textit{primitive} if for any vertices $\la_1,\la_2\in \Ga$ there exists a vertex $\mu\in \Ga$ such that $\mu\geq\la_1,\la_2$.
\end{definition}

\begin{theorem}\cites[Theorem 7 p.143, Corollary p.144 ]{wassermann1981}[Theorem 3.14]{safonkin21}\label{th8}
Let $I$ be an ideal of a primitive graded graph $\Ga$. Strictly positive indecomposable finite and semifinite harmonic functions on $\Ga$ are in a bijective correspondence with the similar functions on $I$. This bijection is defined by the restriction of functions on $\Ga$ to the ideal $I$. The inverse map, which produces an indecomposable finite or semifinite harmonic function on $\Ga$ out of any indecomposable finite or semifinite harmonic function on $I$, is given by 
\begin{equation}
    \varphi(\cdot)\mapsto \lim\limits_{N\to\infty}\summ_{\substack{\mu:\mu\in I\\ |\mu|=N}}\dim(\cdot,\mu)\varphi(\mu).
\end{equation}
\end{theorem}

There is a bijective correspondence $I\leftrightarrow\Ga\backslash I$ between ideals and coideals. Let $J$ be a coideal and $I=\Ga\backslash J$ be the corresponding ideal. Then the following conditions are equivalent:
	\begin{enumerate}[label=\theenumi)]
	    \item if $\left\{\mu \;\middle|\; \la\nearrow\mu\right\}\subset I$, then $\la\in I$
	    
	    \item for any $\la\in J$ there exists a vertex $\mu\in J$ such that $\la\nearrow \mu$.
	\end{enumerate}
	
	\begin{definition}\label{def3}
	An ideal $I$ and the corresponding coideal $J$ are called \textit{saturated}, if they satisfy the conditions above. A saturated ideal $I$ is called \textit{primitive}, if for any saturated ideals $I_1,I_2$ such that $I=I_1\cap I_2$ we have $I=I_1$ or $I=I_2$. A saturated coideal $J$ is called \textit{primitive}, if for any saturated coideals $J_1,J_2$ such that $J=J_1\cup J_2$ we have $J=J_1$ or $J=J_2$.
	\end{definition}
	
	The bijection $I\leftrightarrow\Ga\backslash I$ maps primitive saturated ideals to primitive saturated coideals and vice versa.

 \begin{proposition}\cites{Bratteli1972}[Proposition 2.6]{safonkin21}\label{prop4}
    A saturated coideal $J$ of a graded graph is primitive if and only if it is primitive as a graded graph i.e. for any two vertices $\la_1,\la_2\in J$ we can find a vertex $\mu\in J$ such that $\mu\geq\la_1,\la_2$.
 \end{proposition}

 \begin{proposition}\cites[p. 1720 Lemma 12]{versh_ker_85}[Proposition 3.16]{safonkin21}\label{prop3}
	Let $\Ga$ be a graded graph and $\varphi$ be an indecomposable finite or semifinite harmonic function on it. Then $\supp(\varphi)$ is a primitive saturated coideal.
\end{proposition}

\section{Semifinite harmonic functions on the direct product of graded graphs}\label{sec3}

\begin{definition}
By the \textit{direct product} of graded graphs $\left(\Ga_1,\ka_1\right)$ and $\left(\Ga_2,\ka_2\right)$ we mean the graded graph $\left(\Ga_1\times \Ga_2,\ka_{1}\times \ka_{2}\right)$, where
$$\left(\Ga_1\times\Ga_2\right)_k=\bigsqcup\limits_{\substack{n,m\geq0: \\ n+m=k}}\left(\Ga_1\right)_n\times\left(\Ga_2\right)_m$$
and 
$$\left(\ka_1\times\ka_2\right)\Bigl((\la_1,\mu_1);(\la_2,\mu_2)\Bigr)=
\begin{cases}
    \ka_1(\la_1,\la_2),\ &\text{if}\ \mu_1=\mu_2,\\
    \ka_2(\mu_1,\mu_2),\ &\text{if}\ \la_1=\la_2,\\
    0\ &\text{otherwise.}
\end{cases}$$
\end{definition}

\begin{definition}
    We say that a finite harmonic function $\varphi$ on a branching graph is \textit{normalized} if $\varphi(\diameter)=1$.
\end{definition}

\begin{notation}
We denote by $\exFH(\Ga)$ the set of all normalized indecomposable finite harmonic functions on a branching graph $\left(\Ga,\ka\right)$.
\end{notation}

The next theorem describes the indecomposable finite harmonic functions on the direct product of two branching graphs. Note that the statement is almost obvious for the so-called multiplicative branching graphs. 

The paper \cite{safonkin21} contains a sketch of proof of the following theorem. For reader's convenience we present the argument in full details in Appendix \ref{appendix1}.
 
\begin{theorem}\label{th7}\cite[Proposition A.4]{safonkin21}
    Let $\Ga_1$ and $\Ga_2$ be branching graphs and $\varphi$ be a normalized indecomposable finite harmonic function on $\Ga_1\times\Ga_2$, i.e. $\varphi\in\exFH(\Ga_1\times\Ga_2)$. Then only one of the following situations can occur: 
    \begin{enumerate}[label=\theenumi)]
    \item\label{th7case1} There exist $\varphi_1\in\exFH(\Ga_1)$, $\varphi_2\in\exFH(\Ga_2)$ and real positive numbers $w_1,w_2$ with $w_1+w_2=1$ such that  
\begin{equation}\label{f25}
\varphi(\la,\mu)=w_1^{|\la|}w_2^{|\mu|}\varphi_1(\la)\varphi_2(\mu).   
\end{equation} 
Moreover, these $\varphi_1$, $\varphi_2$, $w_1,w_2$ are uniquely defined.

    \item\label{th7case2} There exist $\varphi_1\in\exFH(\Ga_1)$ such that
    \begin{equation}\label{f26}
        \varphi(\la,\mu)=\begin{cases}
            0,&\text{if}\ \mu\neq\diameter,\\
            \varphi_1(\la),&\text{if}\ \mu=\diameter.
        \end{cases}    \end{equation}

    \item\label{th7case3} There exist $\varphi_2\in\exFH(\Ga_2)$ such that 
    \begin{equation}\label{f27}
        \varphi(\la,\mu)=\begin{cases}
            0,&\text{if}\ \la\neq\diameter,\\
            \varphi_2(\mu),&\text{if}\ \la=\diameter.
        \end{cases}    \end{equation}
    \end{enumerate}
    Furthermore, every harmonic function on $\Ga_1\times\Ga_2$ of the form \ref{th7case1}, \ref{th7case2}, or \ref{th7case3} is indecomposable.
\end{theorem}

\begin{remark}
One can readily see that \eqref{f26} and \eqref{f27} are partial cases of \eqref{f25} corresponding to $w_2=0$ and $w_1=0$. We formulate Theorem \ref{th7} in this form to simplify the comparison with Theorem \ref{th1}.
\end{remark}

Our main goal is to prove a semifinite analog of Theorem \ref{th7}, which deals not only with branching graphs but with arbitrary graded graphs as well, see Theorem \ref{th1}.

\begin{notation}
The set of all indecomposable finite (not identically zero) and semifinite harmonic functions on a graded graph $\Ga$ is denoted by $\exH(\Ga)$. The subset of $\exH(\Ga)$ consisting of strictly positive functions is denoted by $\exHo(\Ga)$. 
\end{notation}

The ideal of a graded graph $\Ga$ generated by a vertex $\la$ will be denoted by $\Ga^{\la}$, i.e. $\Ga^{\la}=\{\mu\in\Ga\mid \mu\geq\la\}$. Such ideals will be called \textit{principal} ideals.

\begin{lemma}\label{lem2}
Let $\Ga_1$ and $\Ga_2$ be primitive graded graphs. If $\varphi\in\exHo(\Ga_1\times\Ga_2)$, then there exist $\varphi_1\in\exHo(\Ga_1)$, $\varphi_2\in\exHo(\Ga_2)$ and unique real positive numbers $w_1,w_2$ with $w_1+w_2=1$ such that
\begin{equation}\label{f16}
\varphi(\la,\mu)=w_1^{|\la|}w_2^{|\mu|}\varphi_1(\la)\varphi_2(\mu).   
\end{equation} 
Moreover, $\varphi_1$ and $\varphi_2$ are defined uniquely up to multiplicative constants.
\end{lemma}
\begin{proof}
Let us show that, if $\varphi$ is given by \eqref{f16}, then numbers $w_1,w_2$ are uniquely defined and $\varphi_1$, $\varphi_2$ are uniquely defined up to multiplicative constants. Assume that 
\begin{equation}\label{f17}
\varphi(\la,\mu)=w_1^{|\la|}w_2^{|\mu|}\varphi_1(\la)\varphi_2(\mu)=\wt{w_1}^{|\la|}\wt{w_2}^{|\mu|}\wt{\varphi_1}(\la)\wt{\varphi_2}(\mu).
\end{equation}

The finiteness ideal of $\varphi_i$ coincides with that of $\wt{\varphi_i}$, $i=1,2$. Then we will assume that all further considerations will be performed inside the finiteness ideal of $\varphi$, which is equal to the direct product of finiteness ideals of $\varphi_1$ and $\varphi_2$. It makes possible to rewrite \eqref{f17} as
\begin{equation}\label{f18}
\left(\cfrac{w_1}{\wt{w_1}}\right)^{|\la|}\cfrac{\varphi_1(\la)}{\wt{\varphi_1}(\la)}=\left(\cfrac{\wt{w_2}}{w_2}\right)^{|\mu|}\cfrac{\wt{\varphi_2}(\mu)}{\varphi_2(\mu)}.
\end{equation}

The left-hand side of \eqref{f18} depends only on $\la$, but not on $\mu$ and the right-hand side depends only on $\mu$, but not on $\la$. Then both sides of \eqref{f18} are constant, hence
\begin{equation}\label{f19}
\wt{\varphi_1}(\la)=c_1\cdot\left(\cfrac{w_1}{\wt{w_1}}\right)^{|\la|}\varphi_1(\la)
\end{equation}
\begin{equation}\label{f20}
\wt{\varphi_2}(\mu)=c_2\cdot\left(\cfrac{w_2}{\wt{w}_2}\right)^{|\mu|}\varphi_2(\mu)
\end{equation}
for some real positive constants $c_1,c_2$.

Next, $\wt{\varphi_1}$ is a harmonic function with respect to the multiplicity function $\varkappa_1$, but on the right-hand side of \eqref{f19} we have a harmonic function with respect to the multiplicity function $\cfrac{w_1}{\wt{w_1}}\cdot\varkappa_1$. Hence $\wt{w}_1=w_1$. Thus, $\varphi_1$ and $\wt{\varphi_1}$ are proportional. Similarly we show that $\wt{w}_2=w_2$ and that $\varphi_2$ and $\wt{\varphi_2}$ are proportional.

Let $\varphi\in\exHo(\Ga_1\times\Ga_2)$ and take any $(\la,\mu)\in\Ga_1\times\Ga_2$ such that $\varphi(\la,\mu)<+\infty$. This pair $(\la,\mu)$ will be fixed till the end of the proof. By Theorem \ref{th7} we have 
\begin{equation}
    \cfrac{\varphi(\nu_1,\nu_2)}{\varphi(\la,\mu)}=w_1^{|\nu_1|-|\la|}w_2^{|\nu_2|-|\mu|}\psi_1(\nu_1)\psi_2(\nu_2)
\end{equation}
for any $\nu_1\geq\la$, $\nu_2\geq\mu$, some numbers $w_1,w_2\in\mathbb{R}_{>0}$ with $w_1+w_2=1$, and some finite strictly positive normalized harmonic functions $\psi_1$ and $\psi_2$ on $\Ga_1^{\la}$ and $\Ga_2^{\mu}$. 

Let us denote by $\varphi_1$ and $\varphi_2$ the extensions of $\psi_1$ and $\psi_2$ to primitive graded graphs $\Ga_1$ and $\Ga_2$ respectively provided by Theorem \ref{th8}. We will show that for any $\nu_1\in\Ga_1$ and any $\nu_2\in\Ga_2$ 
\begin{equation}\label{f08}
    \varphi(\nu_1,\nu_2)\geq\varphi(\la,\mu)\cdot w_1^{|\nu_1|-|\la|}w_2^{|\nu_2|-|\mu|}\varphi_1(\nu_1)\varphi_2(\nu_2),
\end{equation}
Then the claim follows from the indecomposability of $\varphi$.

First, we write
\begin{equation}
    \begin{multlined}
    \varphi(\nu_1,\nu_2)=\limm_{N\to+\infty}\summ_{\substack{\la_1\in\Ga_1,\la_2\in\Ga_2\colon\\ |\la_1|+|\la_2|=N}}\dim\Bigl((\nu_1,\nu_1),(\la_1,\la_2)\Bigr)\varphi(\la_1,\la_2)=\\
    \limm_{N\to+\infty}\summ_{\substack{n_1,n_2\colon\\ n_1+n_2=N}}\ \ \summ_{\substack{\la_1\in\Ga_1,\la_2\in\Ga_2\colon\\|\la_1|=|\nu_1|+n_1,\\|\la_2|=|\nu_2|+n_2}}\binom{N}{n_1}\dim_1(\nu_1,\la_1)\dim_2(\nu_2,\la_2)\varphi(\la_1,\la_2).
    \end{multlined}
    \end{equation}
Omitting all terms except those for which $\la_1\geq\la$ and $\la_2\geq\mu$, we obtain
\begin{equation}
\begin{multlined}
    \varphi(\nu_1,\nu_2)\geq \varphi(\la,\mu)\cdot w_1^{|\nu_1|-|\la|}w_2^{|\nu_2|-|\mu|}\times\\
    \limm_{N\to+\infty}\left[\summ_{\substack{n_1,n_2\colon\\ n_1+n_2=N\\ \phantom{|}}}\binom{N}{n_1}w_1^{n_1}w_2^{n_2}\!\summ_{\substack{\la_1\in\Ga_1^{\la}\colon\\ |\la_1|=|\nu_1|+n_1}}\dim_1(\nu_1,\la_1)\psi_1(\la_1)\summ_{\substack{\la_2\in\Ga_2^{\mu}\colon\\ |\la_2|=|\nu_2|+n_2}}\dim_2(\nu_2,\la_2)\psi_2(\la_2)\right]. 
    \end{multlined}
\end{equation}
Next, we omit all summands except those which satisfy $|n_1-w_1N|<N^{2/3}$. Note that the expressions $\summ_{\substack{\la_1\in\Ga_1^{\la}\colon\\ |\la_1|=|\nu_1|+n_1}}\dim_1(\nu_1,\la_1)\psi_1(\la_1)$ and $\summ_{\substack{\la_2\in\Ga_2^{\mu}\colon\\ |\la_2|=|\nu_2|+n_2}}\dim_2(\nu_2,\la_2)\psi_2(\la_2)$ are non-decreasing in $n_1$ and $n_2$. Hence we can bound $\varphi(\nu_1,\nu_2)$ from below
\begin{equation}
    \begin{multlined}
     \varphi(\nu_1,\nu_2)\geq \varphi(\la,\mu)\cdot w_1^{|\nu_1|-|\la|}w_2^{|\nu_2|-|\mu|}\limm_{N\to+\infty}\left[\summ_{\substack{n_1,n_2\colon\\ n_1+n_2=N\\ |n_1-w_1N|<N^{2/3}}}\binom{N}{n_1}w_1^{n_1}w_2^{n_2}\right.\\
    \left.\summ_{\substack{\la_1\in\Ga_1^{\la}\colon\\ |\la_1|=|\nu_1|+\lfloor w_1N-N^{2/3}\rfloor}}\dim_1(\nu_1,\la_1)\psi_1(\la_1)\summ_{\substack{\la_2\in\Ga_2^{\mu}\colon\\ |\la_2|= |\nu_2|+\lfloor w_2N-N^{2/3}\rfloor}}\dim_2(\nu_2,\la_2)\psi_2(\la_2)\right],
    \end{multlined} 
\end{equation}
where $\lfloor\cdot\rfloor$ is the floor function.

Finally, by the central limit theorem we have
\begin{equation}
    \limm_{N\to+\infty}\summ_{\substack{n_1,n_2\colon\\ n_1+n_2=N\\ |n_1-w_1N|<N^{2/3}}}\binom{N}{n_1}w_1^{n_1}w_2^{n_2}=1.
\end{equation}

Thus, \eqref{f08} follows immediately from the very definition of $\varphi_1$ and $\varphi_2$ and their strict positivity.
\end{proof}

Now we would like to prove Main Theorem \ref{th1}. To do this we will need the following lemmas, see Proposition 5.3 from \cite{safonkin21}, \cite[p.149 Boyer's lemma]{wassermann1981}, and \cite[Theorem 1.10]{boyer83}. 

\begin{lemma}[Boyer's lemma]\label{lem4}
Let $\Ga$ be a graded graph and $\varphi$ be a harmonic function on it. Assume that $I\subset \Ga$ is an ideal, $J=\Ga\backslash I$ is the corresponding coideal and we are given a fixed vertex $\la\in J$. Suppose that there exists a vertex $\la'\in I$ and a positive real number $\beta_{\la}$ such that $\varphi(\la')>0$ and for any vertex $\eta\in I$ lying on a large enough level the following inequality holds
\begin{equation}\label{f06}
\summ_{\mu\in J}\dim(\la,\mu)\varkappa(\mu,\eta)\geq \beta_{\la}\dim(\la',\eta).
\end{equation}

Then $\varphi(\la)=+\infty$. If in addition \eqref{f28} holds for $a=\la'$, then \eqref{f28} holds for $a=\la$ as well.
\end{lemma}

\begin{lemma}\cite[Lemma A.2]{safonkin21}\label{lem1}
Let $\Ga_1$ and $\Ga_2$ be graded graphs. If $J\subset\Ga_1\times\Ga_2$ is a saturated primitive coideal, then it is of the form $J=J_1\times J_2$ for some coideals $J_1\subset\Ga_1$ and $J_2\subset \Ga_2$ such that  
	\begin{itemize}
	\item $J_1,J_2$ are saturated and primitive or
	        
	\item $J_1$ is principle\footnote{We say that a coideal of a graded graph $\Ga$ is principal if it is of the form $\{\la\in\Ga \mid \la\leq\mu\}$ for some $\mu\in\Ga$.} and $J_2$ is saturated and primitive or
	        
	\item $J_1$ is saturated and primitive and $J_2$ is principle. 
	\end{itemize}  
	\end{lemma}
	
\begin{proof} We can form $J_1$ and $J_2$ via natural projections $\Ga_1\times \Ga_2\rightarrow \Ga_1$ and $\Ga_1\times \Ga_2\rightarrow \Ga_2$. Then clearly $J\subset J_1\times J_2$, but primitivity of $J$ implies that $J_1\times J_2\subset J$. Thus, $J=J_1\times J_2$. It is easy to check that $J_1$ and $J_2$ are coideals satisfying the following condition: if $\la,\mu\in J_1$ (or $\in J_2$), then there exists $\nu\in J_1$ (or $\in J_2$) such that $\nu\geq \la,\mu$. By Proposition \ref{prop4} it remains to show that $J_1$ (and similarly $J_2$) is either saturated or principal. Let us check that if $J_1$ is not saturated then it is principal. If $J_1$ is not saturated then there exists a vertex $\la\in J_1$ such that none of the vertices lying above $\la$ belongs to $J_1$. Let $\mu$ be an arbitrary vertex $\mu\in J_1$. Then the property above states that there exists a vertex $\nu\in J_1$ such that $\nu\geq\la,\mu$. But this implies that $\nu=\la$. Thus, $J_1$ is principal.
\end{proof}

\begin{remark}
Lemma \ref{lem1} slightly differs from \cite[Lemma A.2]{safonkin21}. The former deals with arbitrary graded graphs, while the latter deals only with branching graphs. The proof presented here is more straightforward than that given in \cite{safonkin21}. Moreover, we do not use the explicit description of the primitive saturated coideals of a branching graph in terms of the path space, see \cites{strat_voic1975}[p.129]{wassermann1981}[Proposition 2.6]{safonkin21}.
\end{remark}

\begin{proof}[Proof of Theorem \ref{th1}]
Let us check that expressions on the right hand sides of \eqref{f09}, \eqref{f10}, and \eqref{f11} define finite or semifinite harmonic functions on $\Ga_1\times\Ga_2$. Functions from \eqref{f10} and \eqref{f11} are from $\exH(\Ga_1\times\Ga_2)$ by Boyer's lemma, see the proof below. To see that the function from \eqref{f09} is in $\exH(\Ga_1\times\Ga_2)$ note that a harmonic function $\Phi$ is semifinite if and only if it is not finite and the following identity holds (see \cite[Proposition 3.7]{safonkin21})
\begin{equation}\label{f29}
   \Phi(\la)=\lim_{N\to\infty}\summ_{\substack{\mu\colon\mu\geq \la,\ |\mu|=N\\ 0<\varphi(\mu)<+\infty}}\dim(\la,\mu)\Phi(\mu).
\end{equation}

To prove this equality for the function from \eqref{f09} we use an argument with the central limit theorem similar to that which was used in the proof of Lemma \ref{lem2}. Namely, we assume that $\varphi(\la,\mu)=+\infty$ and bound the right hand side of the equality \eqref{f29} from below to see that it equals $+\infty$.

It is clear that the functions defined by \eqref{f10} and \eqref{f11} are indecomposable, if the functions $\varphi_1$ and $\varphi_2$ are. We will check that the function from \eqref{f09} is neccesarily indecomposable, if $\varphi_1$ and $\varphi_2$ are indecomposable. For this we note that $\supp(\varphi)=\supp(\varphi_1)\times\supp(\varphi_2)$ and the direct product of two primitive graded graphs is primitive as well. Thus, by Proposition \ref{prop3} $\supp(\varphi)$ is a primitive graded graph. It is enough to check that $\varphi$ is indecomposable being restricted to $\supp(\varphi)$, hence further considerations will be performed inside $\supp(\varphi)$. From Theorem \ref{th8} and Lemma \ref{lem3} applied to any principal ideal of $\supp(\varphi)$ that lies in the finiteness ideal of $\varphi$, it follows that $\varphi$ is indecomposable being restricted to any principal ideal of $\supp(\varphi)$ that lies in the finiteness ideal. Recall that a principal ideal is an ideal generated by some vertex $\mu$, i.e $\{\la\mid \la\geq \mu\}$. Thus, if $\varphi\geq \psi$, then $\psi=c_I\cdot \varphi$ on any principal ideal $I\subset \supp(\varphi)$ that lies in the finiteness ideal of $\varphi$, where $c_I$ is a positive constant. Finally, we remark that in a primitive graph every two ideals have a non-empty intersection. Thus, $c_I$ does not depend on $I$ and $\varphi$ is indecomposable on $\supp(\varphi)$.

Now we will prove that every harmonic function $\varphi\in\exH(\Ga_1\times\Ga_2)$ is of the form \eqref{f09}, \eqref{f10} or \eqref{f11}. Note that, if $\varphi\in\exH(\Ga_1\times \Ga_2)$, then by Proposition \ref{prop3} and Lemma \ref{lem1} the support $\supp \varphi$ is of the form $J_1\times J_2$, where either $J_1$ and $J_2$ are primitive saturated coideals of $\Ga_1$ and $\Ga_2$ respectively, or one of them is saturated and primitive and another one is principal. The first case of the theorem corresponds to the first case mentioned above and follows from Lemma \ref{lem2} immediately. Let us assume that $J_1\subseteq \Ga_1$ is a saturated primitive coideal and $J_2\subset\Ga_2$ is a principal coideal corresponding to some vertex $\nu_2\in\Ga_2$. Then $\varphi(\la,\mu)=0$ if $\mu\not\leq\nu_2$. Consider the ideal of $J_1\times J_2$ consisting of all pairs $(\la,\nu_2)\in J_1\times J_2$, whose second component equals exactly $\nu_2$. Obviously, it is isomorphic to $J_1$ as a Bratteli diagram. Let us denote this ideal by $\left(J_1\times J_2\right)^{(-,\nu_2)}$. Then Theorem \ref{th8} provides us a bijection between strictly positive finite and semifinite indecomposable harmonic functions on $J_1\times J_2$ and on $\left(J_1\times J_2\right)^{(-,\nu_2)}\simeq J_1$. The last thing which remains to show is to indicate how we should extend harmonic functions from the ideal $\left(J_1\times J_2\right)^{(-,\nu_2)}$ to the whole graph $J_1\times J_2$. If $J_2=\{\nu_2\}$ is a singleton, then the extension is trivial, otherwise it is sufficient to take any $\nu'$ such that $\nu'\nearrow\nu_2$ and prove that $\varphi(\la,\nu')=+\infty$, $\la\in J_1$. Consider the ideal $\left(J_1\times J_2\right)^{(-,\nu')}$ of $J_1\times J_2$ consisting of all pairs $(\la,\mu)\in J_1\times J_2$, whose second component $\mu$ is greater than or equal to $\nu'$, that is $\mu=\nu'$ or $\mu=\nu_2$. Then
\begin{equation}
    \left(J_1\times J_2\right)^{(-,\nu')}=J_1\times\{\nu'\}\sqcup J_1\times\{\nu_2\},
\end{equation}
and $J_1\times\{\nu_2\}$ is an ideal of $\left(J_1\times J_2\right)^{(-,\nu')}$, and there is a natural map $J_1\times\{\nu'\}\rightarrow J_1\times\{\nu_2\}$, $(\la,\nu')\mapsto(\la,\nu_2)$ which, for a trivial reason, is an isomorphism of Bratteli diagrams. Take $\la\in J_1$. Using the notation from Lemma \ref{lem4}, we set $\Ga=\left(J_1\times J_2\right)^{(-,\nu')}$, $I=J_1\times\{\nu_2\}$, $\beta_{(\la,\nu')}=\varkappa_2(\nu',\nu_2)$ and consider $(\la,\nu_2)$ as $\la'$ from the lemma. Then the left hand side of \eqref{f06} consists of only one summand and the inequality turns into the equality, hence Lemma \ref{lem4} implies that $\varphi(\la,\nu')=+\infty$ for any $\nu'\nearrow \nu_2$. By assumption \eqref{f28} holds for $a=(\la,\nu_2)$, hence it holds for $a=(\la,\nu')$ as well.
\end{proof}

\begin{remark}
    One can deal with the cases \ref{th1case2} and \ref{th1case3} from Theorem \ref{th1} applying the Vershik-Kerov ergodic method, see \cite[Theorem on p. 60]{kerov_2003}. This argument turns out to be simplier then that presented above. The reason why we use Boyer's lemma is that it guarantees that the functions defined by the right hand sides of \eqref{f10} and \eqref{f11} are finite or semifinite for any $\varphi_1\in\exH(\Ga_1)$, $\varphi_2\in\exH(\Ga_2)$.
\end{remark}

\subsection{Semifinite Vershik-Kerov ring theorem}

There is a well known analog of the Vershik-Kerov ring theorem for semifinite indecomposable harmonic functions on a multiplicative branching graph, see \cite[Theorem p.134]{vershik_kerov83}, \cite[Proposition 8.4]{gnedin_olsh2006}, and \cite[Theorem p.144]{vershik_kerov83}. We formulate it below (Theorem \ref{th6}) and apply to the case of the direct product of two multiplicative graphs (Corollary \ref{cor3}).

\begin{definition}\cite[p.40]{versh_ker_85}\label{def1}
	A branching graph $\Ga$ is called \textit{multiplicative}, if there exists an associative $\mathbb{Z}_{\geq 0}$-graded $\mathbb{R}$-algebra $A=\bigoplus\limits_{n\geq0} A_n$, $A_0=\mathbb{R}$ with a distinguished basis of homogeneous elements $\{a_{\la}\}_{\la\in\Ga}$, that satisfy the following conditions
	\begin{enumerate}[label=\theenumi)]
	    \item $\deg{a_{\la}}=|\la|$;
	 
	    \item $a_{\diameter}$ is the identity in $A$;
	    
	    \item for $\wh{a}=\summ_{\nu\in\Ga_1}\ka(\diameter,\nu)a_{\nu}$ and any vertex $\la\in\Ga$ we have $\wh{a}\cdot a_{\la}=\summ_{\mu:\la\nearrow\mu}\ka(\la,\mu)a_{\mu}$.
	\end{enumerate}
	
	Moreover, we assume that the structure constants of $A$ with respect to the basis $\{a_{\la}\}_{\la\in\Ga}$ are non-negative. 
\end{definition}

\begin{remark}
    If $(\Ga_1,A,\{a_{\la}\}_{\la\in\Ga_1})$ and $(\Ga_2,B,\{b_{\mu}\}_{\mu\in\Ga_2})$ are multiplicative graphs, then $(\Ga_1\times\Ga_2,A\otimes B,\{a_{\la}\otimes b_{\mu}\}_{\la\in\Ga_1,\mu\in\Ga_2})$ is a multiplicative graph as well, and the element $\wh{a}_{A\otimes B} = \wh{a}_A\otimes\mathds{1}_B+\mathds{1}_A\otimes \wh{a}_B$ plays the role of $\wh{a}$ in Definition \ref{def1}.

\end{remark}

\begin{theorem}[Vershik-Kerov Ring Theorem]\cites[Theorem p.134]{vershik_kerov83}[Proposition 8.4]{gnedin_olsh2006}
A finite normalized harmonic function $\varphi$ on a multiplicative branching graph $\Ga$ is indecomposable if and only if $\varphi(a_{\la}\cdot a_{\mu})=\varphi(a_{\la})\cdot \varphi(a_{\mu})$ for any $\la,\mu\in \Ga$.
\end{theorem}

\begin{theorem}\cite[Theorem p.144]{vershik_kerov83}\label{th6}
For any semifinite indecomposable harmonic function $\varphi$ on a multiplicative branching graph $\Ga$ there exists a finite normalized indecomposable harmonic function $\psi$, such that $\varphi(a_{\la}\cdot a_{\mu})=\psi(a_{\la})\cdot\varphi(a_{\mu})$ for any $\la,\mu\in \Ga$ with $\varphi(a_{\mu})<+\infty$.
\end{theorem}

Now we can associate a \textit{finite} indecomposable harmonic function to any harmonic function $\varphi\in\exH(\Ga)$ on a multiplicative graph $\Ga$. We will denote it by $\varphi^{fin}$, i.e. $\varphi^{fin}=\varphi$, if $\varphi$ is finite and $\varphi^{fin}=\psi$, where $\psi$ is given by Theorem \ref{th6}, if $\varphi$ is semifinite. 

\begin{corollary}\label{cor3}
Let $(\Ga_1,A,\{a_{\la}\}_{\la\in\Ga_1})$ and $(\Ga_2,B,\{b_{\mu}\}_{\mu\in\Ga_2})$ be multiplicative graphs, and $\varphi\in\exH(\Ga_1\times\Ga_2)$. Then 
\begin{itemize}
    \item in case \ref{th1case1} of Theorem \ref{th1}
    \begin{equation}
    \varphi^{fin}(\la,\mu)=w_1^{|\la|}w_2^{|\mu|}\varphi_1^{fin}(\la)\varphi_2^{fin}(\mu);
    \end{equation}
    
    \item in case \ref{th1case2} of Theorem \ref{th1}
    \begin{equation}
        \varphi^{fin}(\la,\mu)=\begin{cases}
            0,&\text{if}\ \mu\neq \diameter,\\
            \varphi_1^{fin}(\la),&\text{if}\ \mu=\diameter;
        \end{cases}
    \end{equation}
    
    \item in case \ref{th1case3} of Theorem \ref{th1}
    \begin{equation}
        \varphi^{fin}(\la,\mu)=\begin{cases}
            0,&\text{if}\ \la\neq \diameter,\\
            \varphi_2^{fin}(\mu),&\text{if}\ \la=\diameter.
        \end{cases}
    \end{equation}
\end{itemize}
\end{corollary}

\section{Slow graphs and inverse symmetric semigroups}\label{sec4}

Let $\NN$ be the representation of natural numbers as a half-line Bratteli diagram, with a single vertex on each level and a single edge between adjacent levels. Recall that for a graded graph~$\Gamma$ the corresponding \emph{slow graph} can be defined as a product $\NN \times \Gamma$. Such graphs were defined and studied in~\cite{Vershik-Nikitin}. Their name comes from the following natural description of the path space for a slow graph. Informally, we take any path in the initial graph, and at each step we either move along the path, or stay at the same vertex once again. More formally, let $\T(\Ga)$ stand for the space of infinite paths in $\Ga$ starting at $\diameter$. We take any path $\left(\diameter\nearrow\la_1\nearrow \la_2\nearrow \ldots\right)\in\T(\Ga)$ and an increasing sequence of positive numbers $i_1 < i_2 < \cdots$, and construct a path 
\begin{multline}
\Bigl((0,\diameter)\nearrow (1,\diameter)\nearrow \ldots
(i_1, \diameter)\nearrow
(i_1, \la_1) \nearrow  (i_1+1, \la_1) \nearrow\ldots\bigr. \\
\Bigl. (i_2, \la_1) \nearrow
(i_2,\la_2) \nearrow \ldots\Bigr)\in\T(\NN \times \Gamma).
\end{multline}

Indecomposable finite harmonic functions are known for the slow graphs.

\begin{theorem}\cite{Vershik-Nikitin}\label{th-slow-fin}
    Let $\Ga$ be a branching graph and $\Phi$ be a normalized indecomposable finite harmonic function on $\NN\times\Ga$, i.e. $\Phi\in\exFH(\NN\times\Ga)$. Then only one of the following two situations can occur: 
    \begin{enumerate}[label=\theenumi)]
        \item\label{l1-th-slow-fin} either there exist $\varphi\in\exFH(\Ga)$ and a real number $w$, $0<w\le 1$, such that  
        \begin{equation*}
            \Phi(n,\mu)=(1-w)^{n}w^{|\mu|}\varphi(\mu)   
        \end{equation*} 
        (and $\varphi$ and $w$ are uniquely defined);
    
        \item\label{l2-th-slow-fin} or $\Phi(n,\mu) = 1$ if $\mu=\diameter$, and $\Phi(n,\mu) = 0$ otherwise.
    \end{enumerate}
    Furthermore, every such harmonic function on $\NN \times\Ga$ is indecomposable.
\end{theorem}
\begin{remark}
    Theorem~\ref{th-slow-fin} easily follows from Theorem~\ref{th7}. The current statement differs slightly from the original one in~\cite{Vershik-Nikitin}~--- the degeneration of the harmonic functions for $w=0$ was not mentioned there. Namely the case~\ref{l2-th-slow-fin} in Theorem~\ref{th-slow-fin} can be obtained from \ref{l1-th-slow-fin} if we set $w=0$, but in this case there is no dependence on $\varphi$. Moreover, the notation in~\cite{Vershik-Nikitin} differs slightly from the notation in the present paper. Namely, there the vertices of the $n$-th level of the corresponding slow graph were denoted by $\{(n,\la) | \la\in\Ga, |\la| \le n\}$ instead of $(\NN\times\Ga)_n = \{(k,\la) | \la\in\Ga, k+|\la|=n\}$.
\end{remark}

Indecomposable semifinite harmonic functions for a slow graph can be easily described as a corollary to Theorem~\ref{th1}. Some simplifications are due to the fact that there is only one (trivial) normalized indecomposable finite harmonic function on $\NN$ and no semifinite ones.

\begin{corollary}\label{cor-slow-semifin}
    Let $\Ga$ be a graded graph and $\Phi$ be an indecomposable semifinite harmonic function on $\NN\times\Ga$, then only one of the following situations can occur: 
    \begin{enumerate}[label=\theenumi)]
    \item There exist an indecomposable semifinite harmonic function $\varphi$ on $\Ga$ and a real number $w$,  $0<w\le1$ such that
        \begin{equation*}
            \Phi(n,\mu)=(1-w)^n w^{|\mu|} \varphi(\mu).  
        \end{equation*} 
        Moreover, $\varphi$ and $w$ are uniquely defined.

    \item There exist $\nu\in\Ga$, $\nu\neq\diameter$ and a real positive number $c$ such that
        \begin{equation*}
            \Phi(n,\mu)=\begin{cases}
                0,         &\text{if} \ \mu\not\leq\nu,\\
                +\infty,   &\text{if} \ \mu<\nu,\\
                c,         &\text{if} \ \mu=\nu.
            \end{cases}    
        \end{equation*}

    \item There exist $m\in\NN$, $m\geq 1$ and $\varphi\in\exH(\Ga)$ such that 
        \begin{equation*}
            \Phi(n,\mu)=\begin{cases}
                0,&\text{if}\ n>m,\\
                +\infty,&\text{if}\ n<m,\\
                \varphi(\mu),&\text{if}\ n=m.
            \end{cases}    
        \end{equation*}
    \end{enumerate} 
    Furthermore, every such harmonic function on $\NN \times\Ga$ is semifinite and indecomposable.
\end{corollary}

A natural example of a slow graph is given by the representation theory of \emph{inverse symmetric semigroups} $R_n$. The semigroup $R_n$ can be defined as the semigroup of partial bijections of the set $\{1, 2, \dots, n\}$, and it naturally contains the symmetric group $S_n$ (the group of bijections of the same set $\{1, 2, \dots, n\}$). Recall that the Bratteli diagram correspodning to the chain of the group algebras $\{\CC[S_n]\}$ is the Young graph~$\YY$ (see~e.g.~\cite{Vershik-Okounkov}). The semigroup algebras $\CC[R_n]$ are semisimple, and the Bratteli diagram correspodning to the chain of the semigroup algebras $\{\CC[R_n]\}$ is given by the slow graph $\NN\times\YY$ (see~\cite{Halverson}). Description of the semifinite harmonic functions in this case was one of the motivations for the present paper.

Recall that indecomposable finite normalized harmonic functions on $\YY$ are given by the celebrated Thoma's theorem~\cite{Thoma, vershik_kerov1980}, and semifinite traces on $\YY$ were described in~\cite[Theorem 3 on p.27]{vershik_kerov1980} 
and~\cite[Theorem 9 on p.150]{wassermann1981}.
\begin{theorem}\label{th-Thoma}\cite{Thoma, vershik_kerov1980}
Every indecomposable finite normalized harmonic function on the Young graph $\mathbb{Y}$ is of the form $\varphi_{\alpha,\beta}$, where $\alpha$ and $\beta$ are sequences of non-decreasing real non-negative numbers $\alpha=(\alpha_1\geq\alpha_2\geq \ldots\geq 0)$ and $\beta=(\beta_1\geq\beta_2\geq\ldots\geq 0)$ subject to $\summ_{i=1}^{\infty}\alpha_i+\summ_{j=1}^{\infty}\beta_j\leq 1$. The function $\varphi_{\alpha,\beta}$ is defined as follows $\varphi_{\alpha,\beta}(\la)=s_{\la}(\alpha\vert\beta)$, where $s_{\la}(\alpha\vert\beta)$ is the image of the Schur function $s_{\la}$ under the map $\sym\rightarrow \mathbb{R}$, $p_1\mapsto 1$, $p_n\mapsto \summ_{i=1}^{\infty}\alpha_i^n+(-1)^{n-1}\summ_{j=1}^{\infty}\beta_j^{n}$ for $n\geq2$.
\end{theorem}

\begin{theorem}\cites[Theorem 3 on p.27]{vershik_kerov1980}[Theorem 9 on p.150]{wassermann1981}\label{th-Young-semifin}
Every indecomposable semifinite harmonic function on the Young graph $\mathbb{Y}$ is proportional to some $\varphi_{\alpha,\beta}^{\nu}$, where $\nu$ is a non-empty Young diagram and $\alpha$ and $\beta$ are tuples of non-decreasing real positive numbers $\alpha=(\alpha_1\geq\alpha_2\geq \ldots\geq\alpha_k>0)$ and $\beta=(\beta_1\geq\beta_2\geq\ldots\geq\beta_l>0)$ subject to $\summ_{i=1}^{k}\alpha_i+\summ_{j=1}^{l}\beta_j=1$. The function $\varphi^{\nu}_{\alpha,\beta}$ is defined as follows
\begin{equation}
\varphi^{\nu}_{\alpha,\beta}(\la)=\begin{cases}
        0, &\text{if}\ \la\notin \mathbb{Y}_{k,l}^{\nu},\\
        +\infty, &\text{if}\ \la\in \mathbb{Y}_{k,l}^{\nu}\ \text{but}\ \la\ \text{does not cover the flange}\ \nu,\\
        \varphi_{\alpha,\beta}(\la^f), &\text{if}\ \la\in \mathbb{Y}_{k,l}^{\nu}\ \text{and}\ \la\ \text{covers the flange}\ \nu,
        \end{cases}
\end{equation}
    where 
    \begin{itemize}
        \item $\mathbb{Y}_{k,l}^{\nu}$ is the coideal of the Young graph formed by all Young diagrams that can be fitted into the infinite hook consisting of $k$ infinite rows and $l$ infinite columns with an added flange of the form $\nu$ to the corner of the hook, see Figure \ref{fig3};

        \item $\varphi_{\alpha,\beta}$ is the finite indecomposable harmonic function on $\mathbb{Y}$ associated to $(\alpha,\beta)$;

        \item $\la^f=\la-\nu$ is the Young diagram $\la$ with the flange $\nu$ removed.
    \end{itemize}
    \begin{figure}[H]
    \centering
    \import{Pictures}{fig1.tex}
    \caption{An example of $\mathbb{Y}_{k,l}^{\nu}$ for $k=4$, $l=3$ and $\nu=(5,5,2,2,1)$.}
    \label{fig3}
    White rows and columns represent infinite rows and columns.
\end{figure}
\end{theorem}

Now we can easily combine Theorem~\ref{th-Thoma} and Theorem~\ref{th-Young-semifin} with Corollary~\ref{cor-slow-semifin} to describe semifinite harmonic functions on $\NN\times\YY$.

\begin{proposition}\label{prop:inv-semigroup_semifin-traces}
Each semifinite indecomposable harmonic function on $\mathbb{N}\times\Ga$ is proportional to $\Phi_{\alpha,\beta,w}^{\nu}$, $\Phi_{\alpha,\beta}^{k}$, $\Phi_{\alpha,\beta}^{\nu,k}$, or $\Phi^{\mu}$, where

\begin{enumerate}[leftmargin=3ex,label=\theenumi)]
    \item $\Phi_{\alpha,\beta,w}^{\nu}$ is defined by 
    \begin{equation}
        \Phi_{\alpha,\beta,w}^{\nu}(n,\la) = (1-w)^{n}w^{|\la|}\varphi_{\alpha,\beta}^{\nu}(\la)
    \end{equation} 
    for some $w\in(0,1]$ and a semifinite indecomposable harmonic function $\varphi_{\alpha,\beta}^{\nu}$ on the Young graph $\mathbb{Y}$ from Theorem \ref{th-Young-semifin};

    \item $\Phi_{\alpha,\beta}^{k}$ and $\Phi_{\alpha,\beta}^{\nu,k}$ are defined by 
    \begin{equation}
        \Phi_{\alpha,\beta}^{k}(n,\la)=\begin{cases}
            0,&\text{if}\ n>k,\\
            +\infty,&\text{if}\ n<k,\\
            \varphi_{\alpha,\beta}(\la),&\text{if}\ n=k,
        \end{cases}\ \ \ \ \ \ \ 
        \Phi_{\alpha,\beta}^{\nu,k}(n,\la)=\begin{cases}
            0,&\text{if}\ n>k,\\
            +\infty,&\text{if}\ n<k,\\
            \varphi_{\alpha,\beta}^{\nu}(\la),&\text{if}\ n=k,
        \end{cases}     \end{equation}
        for some $\varphi_{\alpha,\beta}$ from Theorem \ref{th-Thoma}, $\varphi_{\alpha,\beta}^{\nu}$ from Theorem \ref{th-Young-semifin}, and an integer $k>0$;

    \item $\Phi^{\mu}$ is defined by 
    \begin{equation}
        \Phi^{\mu}(n,\la)=\begin{cases}
            0,&\text{if}\ \la\not\leq\mu,\\
            +\infty,&\text{if}\ \la<\mu,\\
            1,&\text{if}\ \la=\mu,
        \end{cases}    
    \end{equation}
    for some $\mu\in\YY$, $\mu\neq\diameter$.
    \end{enumerate} 
\end{proposition}

\appendix
\section{}\label{appendix}
\subsection{Proof of Theorem \ref{th7}}\label{appendix1}

In this section we give a proof of Theorem \ref{th7}. Note that \cite{safonkin21} contains only a sketch of proof. Theorem \ref{th7} immediately follows from the two lemmas below. 

\begin{lemma}\label{lem5}
    Let $\Ga_1$ and $\Ga_2$ be branching graphs and $\varphi$ be an indecomposable finite normalized harmonic function on their product, i.e. $\varphi\in\exFH(\Ga_1\times\Ga_2)$. Then $\varphi$ is either of the form \ref{th7case1}, or \ref{th7case2}, or \ref{th7case3} from Theorem \ref{th7}.
\end{lemma}
\begin{proof}
Here we reproduce the argument from the proof of Proposition A.4 from \cite{safonkin21}, which envolves the Vershik-Kerov ergodic method, see \cites[p.20, Theorem 2]{versh_ker_81}[Theorem on p.60]{kerov_2003}. 

Recall that by $\T(\Ga)$ we denote the space of infinite paths in a branching graph $\Ga$ starting at $\diameter$. By \cite[p.60, Theorem]{kerov_2003} there exists a path $\tau=\Bigl((\diameter,\diameter),(\la'_1,\mu'_1),(\la'_2,\mu'_2),\ldots\Bigr)\in\T(\Ga_1\times\Ga_2)$ such that for any $\la\in\Ga_1$, $\mu\in\Ga_2$
\begin{equation}\label{f30}
\varphi(\la,\mu)=\lim\limits_{N\to+\infty}\cfrac{\dim\Bigl((\la,\mu),(\la'_N,\mu'_N)\Bigr)}{\dim\Bigl((\la'_N,\mu'_N)\Bigr)}.    
\end{equation}
Next, we can write 
$$\cfrac{\dim\Bigl((\la,\mu),(\la'_N,\mu'_N)\Bigr)}{\dim\Bigl((\la'_N,\mu'_N)\Bigr)}=\cfrac{\Bigl(|\la'_N|\Bigr)^{\downarrow|\la|}\cdot\ \Bigl(|\mu'_N|\Bigr)^{\downarrow|\mu|}}{\Bigl(|\la'_N|+|\mu'_N|\Bigr)^{\downarrow (|\la|+|\mu|)}}\cdot\cfrac{\dim_1\left(\la,\la'_N\right)}{\dim_1\left(\la'_N\right)}\cdot\cfrac{\dim_2\left(\mu,\mu'_N\right)}{\dim_2\left(\mu'_N\right)},$$
where $x^{\downarrow k}=x(x-1)\ldots (x-k+1)$. 

Passing to appropriate subsequences of vertices in $\tau$ we may assume that the following limits exist
$$\lim\limits_{N\to+\infty}\cfrac{\dim_1\left(\la,\la'_N\right)}{\dim_1\left(\la'_N\right)},\ \lim\limits_{N\to+\infty}\cfrac{\dim_2\left(\mu,\mu'_N\right)}{\dim_2\left(\mu'_N\right)},\ \lim\limits_{N\to+\infty}\cfrac{|\la'_N|}{|\la'_N|+|\mu'_N|},\ \lim\limits_{N\to+\infty}\cfrac{|\mu'_N|}{|\la'_N|+|\mu'_N|}.$$ 
Denote them by $\varphi_1(\la)$, $\varphi_2(\mu)$, $w_1$, and $w_2$. Suppose that $w_1$ and $w_2$ are non-zero, this situation corresponds to the case \ref{th7case1} in Theorem \ref{th7}. Then it is easy to check that $\varphi_1$ and $\varphi_2$ are harmonic functions on $\Ga_1$ and $\Ga_2$. They are indecomposable, since $\varphi$ is indecomposable.

Suppose now that $w_1=0$. Then $w_2=1$ and one can check that $\varphi_2$ is a harmonic function on $\Ga_2$ as before. Equation \eqref{f30} turns into \eqref{f27} from the case \ref{th7case3} in Theorem \ref{th7}. By the same argument as before $\varphi_2$ is indecomposable. The case $w_2=0$ can be dealt with similarly; it corresponds to the case \ref{th7case2} in Theorem \ref{th7}. 
\end{proof}

It is clear that the functions defined in cases \ref{th7case2} and \ref{th7case3} in Theorem \ref{th7} are indecomposable, if so are $\varphi_1$ and $\varphi_2$. Below we prove that the function defined in \ref{th7case1} is indecomposable as well, if $\varphi_1$ and $\varphi_2$ are indecomposable.

For any branching graph $\Ga$ we endow $\exFH(\Ga)$ with the pointwise convergence topology, in which it is a metrizable space, since $\Ga$ is countable.

\begin{lemma}\label{lem3}
    Let $\Ga_1$ and $\Ga_2$ be branching graphs, $\varphi_1$ and $\varphi_2$ be some finite normalized harmonic functions on them, and $w_1$ and $w_2$ be some positive real numbers subject to $w_1+w_2=1$. Then the harmonic function $\varphi$ on the graph $\Ga_1\times\Ga_2$ defined by \eqref{f25} is indecomposable if and only if $\varphi_1$ and $\varphi_2$ are indecomposable.
\end{lemma}
\begin{proof}
    If $\varphi$ is indecomposable, then functions $\varphi_1$ and $\varphi_2$ can not be decomposable by the very definition of $\varphi$. Suppose that $\varphi_1$ and $\varphi_2$ are indecomposable. By Choquet's theorem, see \cite{phelps1966} or \cite[Theorem 9.2]{olshanski2003}, there exists a unique probability measure $P$ on the set $\exFH(\Ga_1\times\Ga_2)$, representing $\varphi$ in the following sence. For any $\la_1\in\Ga_1$ and $\la_2\in\Ga_2$ we have
    \begin{equation}\label{f12}
\varphi(\la_1,\la_2)=\int\limits_{\exFH(\Ga_1\times\Ga_2)} \psi(\la_1,\la_2)P(d\psi).
    \end{equation}
    From Lemma \ref{lem5} it follows that indecomposable harmonic functions on $\Ga_1\times\Ga_2$ are of the following three kinds
    \begin{enumerate}[label=\theenumi)]
        \item\label{t1} $\psi(\la_1,\la_2)=u_1^{|\la_1|}u_2^{|\la_2|}\psi_1(\la_1)\psi_2(\la_2)$ for some real positive numbers $u_1,u_2$ with $u_1+u_2=1$ and some $\psi_1\in\exFH(\Ga_1)$, $\psi_2\in\exFH(\Ga_2)$.

        \item\label{t2} $\psi(\la_1,\la_2)=\begin{cases}
                \psi_1(\la_1),\ &\text{if}\ \la_2=\diameter,\\
                0\ &\text{otherwise}
            \end{cases}\ $  for some $\psi_1\in\exFH(\Ga_1)$.

        \item\label{t3} $\psi(\la_1,\la_2)=\begin{cases}
                \psi_2(\la_2),\ &\text{if}\ \la_1=\diameter,\\
                0\ &\text{otherwise}
            \end{cases}\ $  for some $\psi_2\in\exFH(\Ga_2)$.
    \end{enumerate}
    Note that in the first case we can recover these $u_1,u_2,\psi_1$, and $\psi_2$ as follows
    \begin{align}
&u_1=\summ_{|\la_1|=1}\dim(\la_1)\psi(\la_1,\diameter),\!\!\! 
&& u_2 = 1 - u_1,\\
&\psi_1(\la_1)=\cfrac{1}{u_1^{|\la_1|}}\cdot\psi(\la_1,\diameter),\!\!\! 
&&\psi_2(\la_2)=\cfrac{1}{u_2^{|\la_2|}}\cdot\psi(\diameter,\la_2).
    \end{align}
    
The second and third cases above can be considered as a part of the first one with $u_1=1$ and $u_1=0$ respectively. Then we have a natural map from $\exFH(\Ga_1\times \Ga_2)$ to 
$$
M\vcentcolon = 
\quot{[0,1]\times\exFH(\Ga_1)\times\exFH(\Ga_2)}{\sim},
$$ 
where $(u_1,\psi_1,\psi_2)\sim (\wt{u_1}, \wt{\psi_1},\wt{\psi_2})$ if and only if either $u_1=\wt{u_1}=0$, $\psi_2=\wt{\psi_2}$ or $u_1=\wt{u_1}=1$, $\psi_1=\wt{\psi_1}$. It is readily seen that the map $\exFH(\Ga_1\times\Ga_2)\rightarrow M$ is injective and continuous. 
    
Keeping all the previous discussion in mind we multiply \eqref{f12} by $\dim(\la_1)$ and sum over all $\la_1$ liying on a common level $N$ of the graph $\Ga_1$, and use the fact that 
\begin{equation}
\summ_{|\la_1|=N}\dim(\la_1)\psi_1(\la_1)=1.
\end{equation}
Next we do the same thing for $\la_2$. Then we see by de Finetti's theorem, see Theorems 5.1 and 5.2 in \cite{bor_olsh2017}, that the projection of the measure $P$ to the first coordinate is concentrated at the point $w_1$, hence the integration in \eqref{f12} must run over the set of harmonic functions of the first type \ref{t1}. Then the factors $w_1^{|\la_1|}w_2^{|\la_2|}$ and $u_1^{|\la_1|}u_2^{|\la_2|}$ in the resulting expression cancel out, and, with some abuse of notation, we can write 
\begin{equation}\label{f14}
    \varphi_1(\la_1)\varphi_2(\la_2)=\int\limits_{\exFH(\Ga_1)\times\exFH(\Ga_2)} \psi_1(\la_1)\psi_2(\la_2) P(d\psi).
\end{equation}
Taking $\la_2=\diameter$ in \eqref{f14}, we see that the measure $P$ must be concentrated only at $\varphi_1$. Analogously, $P$ is concentrated only at $\varphi_2$. Thus, $P$ is a delta measure on $\exFH(\Ga_1\times\Ga_2)$.
\end{proof}

\begin{notation}
    Let $\exFHo(\Ga)$ denote the set of strictly positive functions from $\exFH(\Ga)$.  
\end{notation}

Next proposition follows immediately from the proof of Lemma \ref{lem3}.

\begin{proposition}\label{prop1}
Consider the space $$\quot{[0,1]\times\exFH(\Ga_1)\times\exFH(\Ga_2)}{\sim},$$ 
where the equivalence relation $\sim$ is defined by $(w,\varphi_1,\varphi_2)\sim (\wt{w},\wt{\varphi_1},\wt{\varphi_2})$ if and only if either $w=\wt{w}=0$, $\varphi_2=\wt{\varphi_2}$ or $w=\wt{w}=1$, $\varphi_1=\wt{\varphi_1}$.

    The following maps    $$\quot{[0,1]\times\exFH(\Ga_1)\times\exFH(\Ga_2)}{\sim}\longrightarrow \exFH(\Ga_1\times \Ga_2)$$
and    $$(0,1)\times\exFHo(\Ga_1)\times\exFHo(\Ga_2)\longrightarrow \exFHo(\Ga_1\times\Ga_2),$$
defined by $(w,\varphi_1,\varphi_2)\mapsto \varphi$, where $\varphi$ is given by \eqref{f25} with $w_1=w$ and $w_2=1-w$, are homeomorphisms. 
\end{proposition}

\subsection{Product of branching graphs revisited}\label{appendix2}
This section is devoted to one simple (and almost elementary) fact about finite harmonic functions on the direct product of two branching graphs, see Proposition \ref{prop2} and Proposition \ref{prop5} below. From now on we consider arbitrary finite harmonic functions on the direct product of two branching graphs, but not only indecomposable ones.

Let $\Ga_1$ and $\Ga_2$ be branching graphs, $\varphi_1$, $\varphi_2$ be finite normalized harmonic functions on them and let $w_1,w_2\in\mathbb{R}_{>0}$ be such that $w_1+w_2=1$. Then 
\begin{equation}\label{f0}
    \varphi(\la_1,\la_2)=w_1^{|\la_1|}w_2^{|\la_2|}\varphi_1(\la_1)\varphi_2(\la_2)
\end{equation}
is a finite normalized harmonic function on $\Ga_1\times \Ga_2$\footnote{Note that for any finite harmonic function $f$ on the Pascal graph $\mathbb{P}$ the function $\varphi$ defined by $\varphi(\la_1,\la_2)=f(|\la_1|,|\la_2|)\varphi_1(\la_1)\varphi_2(\la_2)$ is harmonic on $\Ga_1\times\Ga_2$.}.

\begin{proposition}\label{prop2}
Keeping the aforementioned notation we can recover $\varphi_1$, $\varphi_2$, $w_1$ and $w_2$ from $\varphi$ by the following formulas
\begin{equation}\label{f01}
w_1^{k_1}w_2^{k_2}=\summ_{|\la_1|=k_1,|\la_2|=k_2}\dim_1(\la_1)\dim_2(\la_2)\varphi(\la_1,\la_2),
\end{equation}
\begin{equation}\label{f02}
    \varphi_1(\la_1)=\summ_{n_2=0}^{\infty}\binom{n_2+|\la_1|-1}{n_2}\summ_{|\la_2|=n_2}\dim_2(\la_2)\varphi(\la_1,\la_2),
\end{equation}
\begin{equation}\label{f03}
    \varphi_2(\la_2)=\summ_{n_1=0}^{\infty}\binom{n_1+|\la_2|-1}{n_1}\summ_{|\la_1|=n_1}\dim_1(\la_1)\varphi(\la_1,\la_2).
\end{equation}
\end{proposition}
\begin{proof}
Identities \eqref{f02} and \eqref{f03} reduce to 
\begin{equation}\label{f15}
    \cfrac{1}{(1-y)^{k+1}}=\summ_{n=0}^{\infty}\binom{n+k}{k}y^n
\end{equation}
and \eqref{f01} is obvious.
\end{proof}

Remark that the right hand side of \eqref{f01} defines a harmonic function on the Pascal graph $\Pas$ for any finite normalized harmonic function $\varphi$ on $\Ga_1\times\Ga_2$. Proposition \ref{prop2} shows that the right hand sides of \eqref{f02} and \eqref{f03} define harmonic functions on $\Ga_1$ and $\Ga_2$, if $\varphi$ is of the form \eqref{f0}. In fact, more general claim holds.

\begin{proposition}\label{prop5}
    The functions on $\Ga_1$ and $\Ga_2$ defined by the right hand sides of \eqref{f02} and \eqref{f03} are finite and harmonic for any finite normalized harmonic function $\varphi$ on $\Ga_1\times\Ga_2$.
\end{proposition}
\begin{proof}
The statement is a trivial consequence of Proposition \ref{prop1}, Proposition \ref{prop2}, and Choquet's theorem \cites{phelps1966}[Theorem 9.2]{olshanski2003}.
\end{proof}
\begin{remark}
    One can prove this proposition using only elementary methods and de Finetti's Theorem, see Theorems 5.1 and 5.2 in \cite{bor_olsh2017}. The key observation is that the following expression
    \begin{equation}\label{f07}
    \summ_{|\la_1|=k_1,|\la_2|=k_2}\dim_1(\la_1)\dim_2(\la_2)\varphi(\la_1,\la_2)
    \end{equation}
    defines a harmonic function on the Pascal graph $\Pas$ for any finite harmonic function $\varphi$ on $\Ga_1\times \Ga_2$. Then by Theorem 5.1 from \cite{bor_olsh2017} \eqref{f07} is a mixture of indecomposables. Next, from the inequality 
    \begin{equation}
        \dim_1(\la_1)\summ_{|\la_2|=n_2}\dim_2(\la_2)\varphi(\la_1,\la_2)\leq \summ_{|\mu|=|\la_1|,|\la_2|=n_2}\dim_1(\mu)\dim_2(\la_2)\varphi(\mu,\la_2),
    \end{equation}
    identity \eqref{f15} and the integral representation of \eqref{f07}, it follows that the expression defined by the right hand side of \eqref{f02} is finite and not exceeding $\cfrac{1}{\dim_1(\la_1)}$. Let us denote it by $\pi_1$, i.e.
    \begin{equation}\label{f21}
        \pi_1(\la)\vcentcolon=\summ_{n_2=0}^{\infty}\binom{n_2+|\la|-1}{n_2}\summ_{|\la_2|=n_2}\dim_2(\la_2)\varphi(\la,\la_2)\leq\cfrac{1}{\dim_1(\la)}<+\infty.
    \end{equation}
    
    To establish the harmonicity condition for $\pi_1$ we prove the following identity by induction on $k$
    \begin{equation}\label{f13}
        \begin{multlined}
    \pi_1(\la)-\summ_{\mu\searrow\la}\varkappa_1(\la,\mu)\pi_1(\mu)=
    \summ_{n_2=k}^{\infty}\binom{n_2+|\la|-1-k}{|\la|-1}\summ_{|\la_2|=n_2}\dim_2(\la_2)\varphi(\la,\la_2)-\\
    \summ_{\mu\searrow\la}\summ_{n_2=k}^{\infty}\binom{n_2+|\mu|-1-k}{|\mu|-1}\summ_{|\la_2|=n_2}\dim_2(\la_2)\varphi(\mu,\la_2)
    \end{multlined}
    \end{equation}
    Finally, we recall that 
    \begin{equation}
        \binom{n_2+|\la|-1-k}{|\la|-1}\leq \binom{n_2+|\la|-1}{|\la|-1},
    \end{equation}
    hence each of the summands in \eqref{f13} tends to $0$ as $k$ goes to $+\infty$.
\end{remark}
\begin{remark}
    Take $\varphi$ as in \eqref{f0}, then $\psi_1(\la)\vcentcolon=\varphi(\la,\diameter)$ is a finite harmonic function on a branching graph $(\Ga_1,\wt{\varkappa_1})$ that is similar to $(\Ga_1,\varkappa_1)$ in the sense of Kerov, see Definition in \S4 from \cite{kerov_89}. The multiplicity function of the new graph differs from $\varkappa_1$ by $w_1$, see \eqref{f0}. Below we prove that the same is true for an arbitrary finite harmonic function on $\Ga_1\times\Ga_2$.
    
    We would like to generalize the key observation from the previous remark a bit. The claim is that the following expression defines a harmonic function on the Pascal graph $\Pas$ for any $\nu_1\in\Ga_1$, $\nu_2\in\Ga_2$, and any finite normalized harmonic function $\varphi$ on $\Ga_1\times\Ga_2$
    \begin{equation}
        \Pi^{\nu_1,\nu_2}(k_1,k_2)\vcentcolon=\summ_{\substack{|\la_1|=|\nu_1|+k_1,\\ \!\!|\la_2|=|\nu_2|+k_2}}\dim_1(\nu_1,\la_1)\dim_2(\nu_2,\la_2)\varphi(\la_1,\la_2).
    \end{equation}
    Next by Theorem 5.1 from \cite{bor_olsh2017} it is represented by a probability measure, say $P^{\nu_1,\nu_2}$, on $[0,1]$, that is 
    \begin{equation}
        \Pi^{\nu_1,\nu_2}(k_1,k_2)=\varphi(\nu_1,\nu_2)\int_{[0,1]}w_1^{k_1}w_2^{k_2}P^{\nu_1,\nu_2}(dw).
    \end{equation}

    Then by \eqref{f15} the expression $\pi_1(\la)$ defined by \eqref{f21} equals
    \begin{equation}
        \varphi(\la,\diameter)\int_{[0,1]}w_1^{-|\la|}P^{\la,\diameter}(dw).
    \end{equation}
    Recall that it is finite for all $\la\in\Ga_1$. 
    
    Thus, the harmonicity condition for $\pi_1$ implies that the function $\psi_1(\la)\vcentcolon=\varphi(\la,\diameter)$ is a finite harmonic function on a branching graph that is similar to $(\Ga_1,\varkappa_1)$ in the sense of Kerov. 
\end{remark}

%% file: Pictures/fig1.tex
\tikzset{every picture/.style={line width=0.75pt}} 

\begin{tikzpicture}[x=0.75pt,y=0.75pt,yscale=-1,xscale=1]

\draw  [draw opacity=0] (61,37.67) -- (256,37.67) -- (256,89.67) -- (61,89.67) -- cycle ; \draw   (74,37.67) -- (74,89.67)(87,37.67) -- (87,89.67)(100,37.67) -- (100,89.67)(113,37.67) -- (113,89.67)(126,37.67) -- (126,89.67)(139,37.67) -- (139,89.67)(152,37.67) -- (152,89.67)(165,37.67) -- (165,89.67)(178,37.67) -- (178,89.67)(191,37.67) -- (191,89.67)(204,37.67) -- (204,89.67)(217,37.67) -- (217,89.67)(230,37.67) -- (230,89.67)(243,37.67) -- (243,89.67) ; \draw   (61,50.67) -- (256,50.67)(61,63.67) -- (256,63.67)(61,76.67) -- (256,76.67) ; \draw   (61,37.67) -- (256,37.67) -- (256,89.67) -- (61,89.67) -- cycle ;
\draw  [draw opacity=0] (100,89.67) -- (100,232.67) -- (61,232.67) -- (61,89.67) -- cycle ; \draw   (100,102.67) -- (61,102.67)(100,115.67) -- (61,115.67)(100,128.67) -- (61,128.67)(100,141.67) -- (61,141.67)(100,154.67) -- (61,154.67)(100,167.67) -- (61,167.67)(100,180.67) -- (61,180.67)(100,193.67) -- (61,193.67)(100,206.67) -- (61,206.67)(100,219.67) -- (61,219.67) ; \draw   (87,89.67) -- (87,232.67)(74,89.67) -- (74,232.67) ; \draw   (100,89.67) -- (100,232.67) -- (61,232.67) -- (61,89.67) -- cycle ;
\draw  [draw opacity=0][fill={rgb, 255:red, 155; green, 155; blue, 155 }  ,fill opacity=1 ] (100,89.67) -- (165,89.67) -- (165,115.67) -- (100,115.67) -- cycle ; \draw   (113,89.67) -- (113,115.67)(126,89.67) -- (126,115.67)(139,89.67) -- (139,115.67)(152,89.67) -- (152,115.67) ; \draw   (100,102.67) -- (165,102.67) ; \draw   (100,89.67) -- (165,89.67) -- (165,115.67) -- (100,115.67) -- cycle ;
\draw  [draw opacity=0][fill={rgb, 255:red, 155; green, 155; blue, 155 }  ,fill opacity=1 ] (100,115.67) -- (126,115.67) -- (126,141.67) -- (100,141.67) -- cycle ; \draw   (113,115.67) -- (113,141.67) ; \draw   (100,128.67) -- (126,128.67) ; \draw   (100,115.67) -- (126,115.67) -- (126,141.67) -- (100,141.67) -- cycle ;
\draw   (255.83,89.71) .. controls (260.5,89.76) and (262.85,87.45) .. (262.9,82.78) -- (262.99,73.73) .. controls (263.05,67.06) and (265.41,63.75) .. (270.08,63.8) .. controls (265.41,63.75) and (263.11,60.4) .. (263.18,53.73)(263.15,56.73) -- (263.27,44.69) .. controls (263.31,40.02) and (261,37.67) .. (256.33,37.62) ;
\draw   (61,232.79) .. controls (61.05,237.46) and (63.4,239.77) .. (68.07,239.72) -- (70.58,239.69) .. controls (77.24,239.62) and (80.6,241.91) .. (80.65,246.58) .. controls (80.6,241.91) and (83.9,239.55) .. (90.57,239.48)(87.57,239.51) -- (93.07,239.45) .. controls (97.74,239.4) and (100.05,237.05) .. (100,232.38) ;
\draw  [draw opacity=0][fill={rgb, 255:red, 155; green, 155; blue, 155 }  ,fill opacity=1 ] (100,141.67) -- (113,141.67) -- (113,154.67) -- (100,154.67) -- cycle ; \draw    ; \draw    ; \draw   (100,141.67) -- (113,141.67) -- (113,154.67) -- (100,154.67) -- cycle ;

\draw (150,128.4) node [anchor=north west][inner sep=0.75pt]    {$\nu $};
\draw (278.69,56.21) node [anchor=north west][inner sep=0.75pt]    {$k$};
\draw (77.6,253.43) node [anchor=north west][inner sep=0.75pt]    {$l$};

\end{tikzpicture}